\documentclass[11pt,leqno]{article}
\usepackage[latin1]{inputenc}
\usepackage{graphicx}
\usepackage{amsfonts}
\usepackage[french]{babel}
\usepackage[T1]{fontenc}
\usepackage{amssymb}
\usepackage{a4wide}

\newenvironment{preuve}
      {\par\noindent{\it Démonstration\/: }\nopagebreak\normalsize}%
				  {\linebreak[2]\hspace*{\fill}
				  $\Box$
				  \ifdim\lastskip<12pt
       \removelastskip \penalty-200  \vskip12pt  \fi}

\newtheorem{thm}{Th\'eorème}[section]
\newtheorem{prop}[thm]{Proposition}

\newtheorem{lem}[thm]{Lemme}

\newtheorem{rem}[thm]{Remarque}

\font\mathsymb=msbm10
\def\Bbb#1{\hbox{\mathsymb #1}}
\def\R{{\Bbb R}}

\def\N{{\Bbb N}}
\def\Q{{\Bbb Q}}

\def\fr{\displaystyle \frac}

\def\vect{\,{\rm Vect}\,}

\def\ib{\begin{itemize}}
\def\ie{\end{itemize}}

\def \N{{\Bbb N}}

\def \fr{\frac}

\def \com{{\mbox{\rm Com}}}

\def \fr{\frac}

\def \air{{\vskip 12pt\noindent}\par}

\def \ib{\begin{enumerate}}
\def \ie{\end{enumerate}}

\def \Id{{\;\mbox{\rm Id}}}

\def \det{{\;\mbox{\rm det}}}
\def \pgcd{{\;\mbox{\rm pgcd}}}

\def \RR {{\Bbb R}}
\def \SS {{\Bbb S}}

\def\ib {\begin{enumerate}}
\def\ie {\end{enumerate}}

\def\adots{\mathinner{\mkern2mu\raise 1pt\hbox{.}\mkern 3mu\raise
4pt\hbox{.}\mkern1mu\raise 7pt\hbox{{.}}}}

\begin{document}
\title{Représentations déterminantales effectives des polynômes univariés par les matrices flèches}
\author{Ronan Quarez\\ \small IRMAR (CNRS, URA 305), Universit\'e de Rennes 1, Campus de Beaulieu\\ 35042 Rennes Cedex, France\\ {\small e-mail : ronan.quarez@univ-rennes1.fr}}

\maketitle

\air\par
{\bf Résumé} : Nous établissons d'abord l'existence d'une représentation déterminantale effective pour tout polynôme univarié à coefficient réels. Nous voyons ensuite plus précisément qu'un polynôme univarié de degré $r+2s$ à coefficients réels admet une représentation déterminantale effective de signature $(r+s,s)$ si et seulement s'il possède au moins $r$ racines réelles comptées avec multiplicité.  Les représentations déterminantales effectives sont construites à partir de matrices flèches. 
\air

{\small {\bf Mots clé} : Corps réel clos - Forme bilinéaire symétrique - Matrice flèche - Représentation déterminantale - Séries de Puiseux}\par {\small {\bf MSC classification} : 12 - 13 - 15}



\section{Introduction}
D'après \cite{HMV}, on peut affirmer que tout polynôme $p(x)\in\R[x_1,\ldots,x_N]$ tel que $p(0)\not =0$, possède une représentation déterminantale : $$p(x)=\alpha\det(J-\sum_{i=1}^Nx_i A_i)$$ où $\alpha\in\RR$, $J$ est une {\it matrice signature} (une matrice diagonale avec des coefficients $\pm 1$ sur la diagonale) de $\RR^{D\times D}$, et $A_1,\ldots A_N$ des matrices symétriques de $\RR^{D\times D}$ (dont on notera l'ensemble $\SS\RR^{D\times D}$). On a bien sûr $\alpha=\det(J)p(0)$ et dans la suite on se ramènera souvent au cas d'un polynôme $p(x)$ normalisé, i.e. tel que $\alpha=1$. \par
\air
Lorsque $N=1$ ou $N=2$ on peut prendre $D=d$ le degré du polynôme $p(x)$ (confer \cite{HV}). Par contre, lorsque $N\geq 3$, on a en général $D>d$. L'existence d'une représentation déterminantale générale étant due à \cite[Theorem 14.1]{HMV}, on peut aussi regarder \cite{Qz} pour une démonstration n'utilisant que de l'algèbre linéaire.
\air
Rappelons de \cite{HV} qu'un polynôme $p(x)$ est dit {\bf RZ}, si la condition $p(\lambda x)=0$ implique que $\lambda$ est réel.
Notons d'abord que tout polynôme admettant une représentation déterminantale {\it unitaire}, c'est-à-dire avec $J=\Id_d$ la matrice identité de $\RR^{d\times d}$ (que l'on note encore $\Id$ lorsqu'il n'y a pas d'ambiguïté) est {\bf RZ}. En effet, 
$\det(\Id-\lambda(\sum_{i=1}^Nx_iA_i))=0$ implique que $\fr{1}{\lambda}$ est une valeur propre de la matrice symétrique  $\sum_{i=1}^Nx_iA_i$, donc $\lambda$ est réel.\par
On peut alors légitimement s'interroger sur la propriété réciproque, à savoir si $p(x)$ est un polynôme {\bf RZ}, est-ce qu'il possède une représentation déterminantale unitaire ?
\par
Si $N\geq 3$, on ne sait pas répondre à la question.
Si $N=2$, le théorème profond \cite[Theorem 2.2]{HV} y répond par l'affirmative.\par
Si $N=1$ le résultat est évident. Il suffit en effet de partir  de la décomposition du polynôme $p(x)$ en produit d'irréductibles. Comme $p(x)$ est {\bf RZ}, il n'a que des racines réelles, donc il est de la forme $p(x)= \prod_{i=1}^d(1-\alpha_ix)$ pour peu qu'on suppose $p(x)=1$. On en tire l'identité $$p(x)=\det\left(\Id-x\left(\begin{array}{ccc}
\alpha_1&0&0\\
0&\ddots&0\\
0&0&\alpha_d
\end{array}\right)\right)$$ 

Toutefois, le défaut de cette représentation est de faire appel aux racines de $p(x)$. Dans la suite, nous nous intéressons aux représentations déterminantales {\it effectives}, c'est-à-dire que l'on peut construire en temps polynomial en fonction de la taille des coefficients et du degré de $p(x)$. Typiquement, nous exhibons des matrices dont les coefficients sont donnés par des formules polynomiales explicites en fonction des coefficients de $p(x)$ ou de réels qui s'en déduisent par des algorithme polynomiaux classiques.

Par commodité, on pourra manipuler le problème équivalent de la recherche d'une matrice symétrique dont le polynôme caractéristique est donné. Il suffit en effet de noter que lorsqu'on fixe la taille de la matrice $A$ égale au degré $d$ du polynôme, la condition $p(x)=\det(\Id-xA)$ équivaut à $p^*(x)=\det(x\Id-A)$ où $p^*(x)$ est le polynôme réciproque de $p(x)$.
Cette question a reçu les contributions de \cite{Fi} où il est fait usage de matrices flèches, et de \cite{Sr} où il est fait usage de matrices tridiagonales. 
\air
Dans ce papier nous reprenons la méthode de \cite{Fi} et nous établissons en \ref{effective_rep_determ_fleche} et \ref{recip} qu'un polynôme $p(x)$ univarié de degré $d=r+2s$ possède au moins $r$ racines réelles comptées avec multiplicité si et seulement s'il admet (sous la forme d'une matrice flèche) une {\it représentation déterminantale effective de taille $d$ et de signature $(r+s,s)$}, c'est-à-dire qu'on peut déterminer de manière effective une matrice symétrique $A\in\RR^{d\times d}$ telle que, à un facteur multiplicatif près, $p(x)=\det(J_s-xA)$, avec $J_s=\Id_{s+r}\bigoplus(-\Id_s)$ où $\bigoplus$ désigne la somme directe usuelle entre matrices.
\air Le plan est le suivant. Au paragraphe 2, nous exposons deux Lemmes techniques concernant le déterminant d'une matrice flèche, ainsi que quelques notions utiles sur les corps réels clos et les séries de Puiseux. Au paragraphe 3, nous rappelons la méthode de \cite{Fi}, méthode que nous généralisons au cas des polynômes quelconques au paragraphe 4.
Puis, au paragraphe 5, nous nous intéressons à la réciproque de cette propriété et nous montrons au passage que tout endomorphisme auto-adjoint relativement à une forme quadratique de signature $(r+s,s)$ possède au moins $r$ valeurs propres réelles comptées avec multiplicité. Pour terminer, au paragraphe 6, nous mentionnons que l'on peut passer  de manière effective d'une représentation déterminantale  de signature $(r+s,s)$ donnée par une matrice flèche à une représentation déterminantale de signature $(r+s-1,s+1)$. 


\section{Déterminant d'une matrice flèche et corps réels clos}
Nous dirons qu'une matrice $A$ donnée par blocs par 
$$A=\left(\begin{array}{cc}
B&H\\
K&D
\end{array}\right)$$ 
 avec $B\in\RR^{r\times r}$, $H\in\RR^{r\times s}$, $K\in\RR^{s\times r}$, $D\in\RR^{s\times s}$ est une {\it matrice flèche} lorsque $B$ est diagonale, ou par extension une matrice diagonale par blocs avec des blocs de taille $1$ ou $2$.
   
\par
Le point clé des calculs matriciels à suivre sur les matrices flèches se fera souvent en considérant leurs {\it compléments de Schur}, ce qu'on explicite dans le Lemme classique suivant :
 
\begin{lem}\label{schur} 
Si $B\in\RR^{r\times r}$ est inversible, $H\in\RR^{r\times s}$, $K\in\RR^{s\times r}$, $D\in\RR^{s\times s}$, alors :
\ib
\item[(i)]   
$$\det\left(\begin{array}{cc}
B&H\\
K&D
\end{array}\right)=\det(B)\det(D-KB^{-1}H)$$ 

En particulier, 
\item[(ii)] lorsque $s=1$, on obtient : 
$$\begin{array}{ccc}
\det\left(\begin{array}{cc}
B&H\\
K&D
\end{array}\right)&=&\det(D)\det(B)-\det(D)KB^{-1}H\\
&=&\det(D)\det(B)-K(\com(B))^TH
\end{array}$$ 
où $\com(B)$ est la comatrice de $B$.
\ie
\air
Lorsque $s=1$ et $B$ est diagonale, on a peut expliciter le polynôme caractéristique suivant :\air  
\ib
\item[(iii)]
$${\footnotesize
\begin{array}{lcc}
\det\left(\begin{array}{ccccc}
x-\lambda_1&0&\ldots&0&h_1\\
0&x-\lambda_2&\ddots&\vdots&h_2\\
\vdots&\ddots&\ddots&0&\vdots\\
0&\ldots&0&x-\lambda_r&h_r\\
k_1&k_2&\ldots&k_r&x-d\\
\end{array}\right)&=&(x-d)\times\prod_{i=1}^r(x-\lambda_i)-\sum_{i=1}^rh_ik_i\prod_{j\not= i}(x-\lambda_j)
\end{array}}$$ 
\ie
\end{lem}
\air

Par la suite, nous seront amené à changer le corps de base $\RR$ en un sur-corps avec qui il partagera toutes les propriétés de corps ordonné. 
Un corps $R$ est dit {\it réel clos} si $R[I]$ est algébriquement clos, où $I$ est tel que $I^2=-1$.\par
Si $R$ est réel clos, alors $R\langle\langle\epsilon\rangle\rangle=\{\sum_{i=-i_0}^{+\infty}a_i\epsilon^{i/q}\mid i_0\in \N,q\in\N^*,\forall i\; a_i\in R\}$ est encore réel clos. C'est le corps des {\it séries de Puiseux} à coefficients dans $R$ en la variable $\epsilon$. Sa clôture algébrique est $R\langle\langle\epsilon\rangle\rangle[I]$ ou encore $R[I]\langle\langle\epsilon\rangle\rangle$. Le plus petit entier relatif $p$ tel que $a_{p}\not =0$ est appelé la {\it valuation} de la série de Puiseux. \par
Si $R$ est réel clos, les polynômes de $R[x]$ vérifient, entre autres, le Théorème des valeurs intermédiaires, le Théorème de Rolle, ainsi que le Théorème de Sturm. 
\par
Nous serons amené aussi à considérer $R\langle\epsilon\rangle$ le sous corps de $R\langle\langle\epsilon\rangle\rangle$ dont chaque élément est algébrique au-dessus de $R(\epsilon)$. Sur le sous-anneau $\RR\langle\epsilon\rangle_0$ constitué des éléments de $\RR\langle\epsilon\rangle$ qui sont de valuation positive (typiquement, si $S(\epsilon)$ est une série de Puiseux qui annule un polynôme {\it unitaire} à coefficients dans $\RR(\epsilon)$, alors elle est de valuation positive), on dispose du morphisme d'évaluation à valeurs dans $R$ qui consiste à substituer $0$ à $\epsilon$. On notera ce morphisme $\lim_{\epsilon\to 0}$.\par 
Pour toutes ces notions, on peut se référer à \cite{BPR}. 
\air
Dans le Lemme qui suit, on peut voir, d'un point de vue empirique,  l'ajout de la nouvelle variable $\epsilon$ comme une perturbation infinitésimale d'une matrice qui permet de se "débarrasser" de ses valeurs propres multiples. \air

\begin{lem}\label{approx_schur} 
Soit $R$ un corps réel clos.
Soit la matrice
$A=\left(
\begin{array}{cc}
B&H\\
K&d
\end{array}\right)$
 où $B\in R^{n\times n}$, $H\in R^{n\times 1}$, $K\in R^{1\times n}$, $d\in R$.
On suppose que $B$ ne possède que des valeurs propres simples dans $R[I]$. Alors, $A_\epsilon=\left(
\begin{array}{cc}
B&H\\
K&d+\epsilon
\end{array}\right)\in R\langle\langle\epsilon\rangle\rangle^{(n+1)\times(n+1)}$
ne possède que des valeurs propres simples dans $R[I]\langle\langle\epsilon\rangle\rangle$. Et on a bien sûr aussi $\lim_{\epsilon\to 0}A_\epsilon=A$.
\end{lem}
\begin{preuve}
Notons $p_A(x)=\det(x\Id-A)$ le polynôme caractéristique de $A$ et $\lambda_1,\ldots,\lambda_{n+1}$ ses racines dans $R[I]$ (i.e. les valeurs propres de $A$). De même, notons $\lambda_{1,\epsilon},\ldots,\lambda_{n+1,\epsilon}$ les valeurs propres de $A_\epsilon$, ou encore les racines de $p_{A_\epsilon}(x)$ dans $R\langle\langle\epsilon\rangle\rangle[I]$. Du fait que $p_{A_\epsilon}(x)$ est unitaire, les $\lambda_{i,\epsilon}$ sont dans $R\langle\epsilon\rangle_0$  pour tout $i=1\ldots n+1$.\par
Comme $\lim_{\epsilon\to 0}p_{A_\epsilon}(x)=p_A(x)$, on sait que $\lim_{\epsilon\to 0}\lambda_{i,\epsilon}$ est une valeur propre de $p_A(x)$, donc à renumérotation près, on peut supposer que $\lim_{\epsilon\to 0}\lambda_{i,\epsilon}=\lambda_i$, pour tout $i=1\ldots n+1$. En particulier, si $\lambda_i$ est une valeur propre simple de $A$, alors $\lambda_{i,\epsilon}$ est une valeur propre simple de $A_{\epsilon}$.

Du Lemme \ref{schur} (ii), on tire

$$p_A(x)=p_B(x)(x-d)-q_H(x)\quad{\rm avec}\quad 
q_H(x)=K(\com(x\Id-B))^TH$$
Et aussi 
$$p_{A_\epsilon}(x)=p_B(x)(x-d-\epsilon)-q_H(x)=p_A(x)-\epsilon p_B(x)$$ \air

Soit $\lambda_\epsilon$ une racine de $p_{A_\epsilon}(x)$ telle que $\lim_{\epsilon\to 0}\lambda_\epsilon=\lambda$ est une racine de $p_A(x)$ de multiplicité $m\geq 2$. \par
Ecrivons $$p_A(x)=\sum_{k=m}^{n+1}\frac{p_A^{(k)}(\lambda)}{k!}(x-\lambda)^k$$ et $$p_B(x)=\sum_{l=0}^{n}\frac{p_B^{(l)}(\lambda)}{l!}(x-\lambda)^l$$
On peut supposer que $p_{A_\epsilon}(\lambda_\epsilon)=p'_{A_\epsilon}(\lambda_\epsilon)=0$ sinon $\lambda_\epsilon$ serait racine simple de
$p_{A_\epsilon}(x)$ et on aurait rien à montrer.
\par

Par ailleurs, on peut écrire $\lambda_\epsilon=\lambda+\alpha\epsilon^u+\ldots$ où $\alpha\in R[I]^*$, $u\in\Q_+^*$, et les points de suspension désignent des termes de valuation $>u$. En effet, on remarque pour cela que $\alpha=0$ est impossible. Sinon, $\lambda_\epsilon=\lambda$ donnerait $p_{A_\epsilon}(\lambda)=-\epsilon p_B(\lambda)=0$ et de même $p'_B(\lambda)=0$, donc $p_B(\lambda)=p'_B(\lambda)=0$ ce qui est contraire aux hypothèses.\par  

Deux cas se présentent :\air
\ib
\item[1)]
Supposons que $p_B(\lambda)\not=0$. Alors 
$$p_{A_\epsilon}(\lambda_\epsilon)=p_{A}(\lambda_\epsilon)-\epsilon p_{B}(\lambda_\epsilon)=\left(\frac{p_A^{(m)}(\lambda)}{m!}\alpha^m\epsilon^{mu}+\ldots\right)-\left(p_B(\lambda)\epsilon+\ldots\right)$$
Donc une condition nécessaire pour que 
$p_{A_\epsilon}(\lambda_\epsilon)=0$ est 
$$mu=1 \quad{\rm et}\quad \frac{p_A^{(m)}(\lambda)}{m!}\alpha^m-p_B(\lambda)=0$$
Puis, on dit que $p_B(x)$ est de la forme $p_B(x)=p_B(\lambda)+\frac{P_B^{(s)}(\lambda)}{s!}x^s+\ldots$ et donc que
$$p'_{A_\epsilon}(\lambda_\epsilon)=p'_{A}(\lambda_\epsilon)-\epsilon p'_{B}(\lambda_\epsilon)=\left(\frac{p_A^{(m)}(\lambda)}{(m-1)!}\alpha^{m-1}\epsilon^{(m-1)u}+\ldots\right)-\left(\frac{p_B^{(s)}(\lambda)}{(s-1)!}\alpha^{s-1}\epsilon^{s-1+1}+\ldots\right)$$
D'où une condition nécessaire pour que $p'_{A_\epsilon}(\lambda_\epsilon)=0$ est 
$$(m-1)u=s\geq 1 \quad{\rm et}\quad \frac{p_A^{(m)}(\lambda)}{(m-1)!}\alpha^{m-1}-\frac{p_B^{(s)}(\lambda)}{(s-1)!}\alpha^{s-1}=0$$
On obtient donc $mu=1$ et $(m-1)u\geq 1$, une contradiction.

\item[2)] Supposons que $p_B(\lambda)=0$.\par Par hypothèse, on a $p'_B(\lambda)\not=0$. L'égalité 
$p_{A_\epsilon}(\lambda_\epsilon)=0$ donne la condition nécessaire :
$$mu=u+1\quad{\rm et}\quad \frac{p_A^{(m)}(\lambda)}{m!}\alpha^{m}=p'_B(\lambda)\alpha$$

De même, l'égalité 
$p'_{A_\epsilon}(\lambda_\epsilon)=0$ donne la condition nécessaire :
$$(m-1)u=1 \quad{\rm et}\quad \frac{p_A^{(m)}(\lambda)}{(m-1)!}\alpha^{m-1}=p'_B(\lambda)$$
D'où $\frac{p_A^{(m)}(\lambda)}{m!}=\frac{p_A^{(m)}(\lambda)}{(m-1)!}$, une contradiction.
\ie\air
Par conséquent, même si $\lambda$ est racine au moins double de $p_A(x)$, alors $\lambda_\epsilon$ n'est plus que racine simple de $p_{A_\epsilon}(x)$. Ainsi, $p_{A_\epsilon}(x)$ ne possède que des racines simples dans $R[I]\langle\langle\epsilon\rangle\rangle$.
\end{preuve}

\section{Représentations déterminantales des polynômes dont toutes les racines sont réelles}
La proposition suivante est issue de \cite{Fi}.
\begin{prop}\label{Fiedler_fleche}
Si $p(x)$ est un polynôme unitaire de degré $d$ scindé à racines simples dans $\RR[x]$, on peut trouver de manière effective une matrice flèche $A\in\RR^{d\times d}$ de la forme $$A=
\left(\begin{array}{ccccc}
\lambda_1&0&\ldots&0&h_1\\
0&\ddots&\ddots&\vdots&\vdots\\
\vdots&\ddots&\ddots&0&\vdots\\
0&\ldots&0&\lambda_{d-1}&h_{d-1}\\
h_1&\ldots&\ldots&h_{d-1}&e
\end{array}\right)
$$
telle que $\det(x\Id-A)=p(x).$ 
\end{prop}

\begin{rem}
Notons déjà que l'on peut étendre le résultat au cas des polynômes ayant toutes leurs racines réelles mais pas forcément simples. Il suffit en effet de factoriser par $\pgcd(p(x),p'(x))$ et de raisonner par récurrence sur le degré $d$ de $p(x)$. Il faut juste prendre garde que l'on obtient au bout du compte une matrice $A$ qui est {\rm flèche-diagonale}, c'est-à-dire diagonale par bloc, chaque bloc étant une matrice flèche.\par
\end{rem}

Rappelons la méthode utilisée dans \cite{Fi}.\par

\air{\bf Localisation des racines}\par
Notons $\alpha_1<\alpha_2<\ldots<\alpha_d$ les racines de $p(x)$. La première étape consiste à localiser les racines, c'est possible par exemple en faisant appel aux suites de Sturm, ou par d'autres moyens équivalents (confer \cite{BPR}). 

\air{\bf Entrelacement}\par
On choisit ensuite  une famille de réels {\it arbitraires} $\lambda_1,\ldots,\lambda_{d-1}$ qui entrelacent les racines, c'est-à-dire tels que  
$$\alpha_1<\lambda_1<\alpha_2<\lambda_2<\ldots<\alpha_{d-1}<\lambda_{d-1}<\alpha_d$$\air   
On cherche alors la matrice $A$ sous la forme :
$$A=
\left(\begin{array}{ccccc}
\lambda_1&0&\ldots&0&h_1\\
0&\ddots&\ddots&\vdots&\vdots\\
\vdots&\ddots&\ddots&0&\vdots\\
0&\ldots&0&\lambda_{d-1}&h_{d-1}\\
h_1&\ldots&\ldots&h_{d-1}&e
\end{array}\right)
$$ 
où $h_1,\ldots,h_{d-1}$ et $e$ sont des réels à déterminer.
 
\air{\bf Interpolation}\par 
De la formule de \ref{schur} (iii), on tire :
$$q(x)=\det(x\Id-A)=(x-e)\prod_{i=1}^{d-1}(x-\lambda_i)-\sum_{i=1}^{d-1}h_i^2\prod_{j\not =i}(x-\lambda_j)$$
Il suffit alors d'interpoler en les réels $\lambda_i$. \par En effet, on a $q(\lambda_i)=-h_i^2\prod_{j\not =i}(\lambda_i-\lambda_j)$, et comme $p(x)$ est unitaire, le réel $p(\lambda_i)$ est du signe de $-\prod_{j\not =i}(\lambda_i-\lambda_j)$ en raison de la condition d'entrelacement. Donc, on peut trouver $h_i$ tel que $q(\lambda_i)=p(\lambda_i)$ pour $i=1\ldots d-1$. De plus, si on choisit $e$ tel que $-e-\sum_{i=1}^{d-1}{\lambda_i}$ vaut le coefficient en $x^{d-1}$ de $p(x)$, on est assuré que $p(x)=q(x)$.


\section{Représentations déterminantales des polynômes quelconques}
Considérons désormais un polynôme $p(x)$ de degré $d$ possédant exactement $r$ racines réelles comptées avec multiplicité.

Encore une fois, plutôt que la recherche d'une représentation déterminantale proprement dite, on traitera le problème équivalent de la recherche d'une matrice symétrique $A$ telle que $p(x)=\det(xJ-A)$ avec $J$ une matrice de signature. 

Pour peu qu'on suppose le coefficient dominant de $p(x)$ égal à $(-1)^s$, on peut l'écrire : $$p(x)=\prod_{i=1}^r(x-\alpha_i)\prod_{j=1}^s(-(x-\beta_j)^2-\gamma_j^2)$$ où les $\alpha_i$, $\beta_j,\gamma_j$ sont réels et $\gamma_j\not =0$.\par

On dispose alors de la représentation évidente 
$p(x)=\det(Jx-A)$
avec 
$$J=\left(\begin{array}{cccccccc}
1&0&&&&&&\\
0&-1&&&&&&\\
&&\ddots&&&&&\\
&&&1&0&&&\\
&&&0&-1&&&\\
&&&&&1&&\\
&&&&&&\ddots\\
&&&&&&&1\\
\end{array}\right)=\left(\Id_1\bigoplus(-\Id_1)\right)^s\bigoplus\Id_r$$
 
et     
$$A=
\left(\begin{array}{cccccccc}
\beta_1&\gamma_1&&&&&&\\
\gamma_1&-\beta_1&&&&&&\\
&&\ddots&&&&&\\
&&&\beta_s&\gamma_s&&&\\
&&&\gamma_s&-\beta_s&&&\\
&&&&&\alpha_1&&\\
&&&&&&\ddots&\\
&&&&&&&\alpha_r\\
\end{array}\right)$$ 

Voyons maintenant comment étendre la méthode du paragraphe précédent pour obtenir une représentation effective.

\begin{thm} \label{effective_rep_determ_fleche}
Tout polynôme $p(x)$ de degré $d=r+2s$ tel que $p(0)\not =0$ et possédant exactement $r$ racines réelles, comptées avec multiplicité, admet une représentation déterminantale effective de signature $(r+s,s)$.
\end{thm}

\begin{preuve}
Encore une fois, quitte à factoriser successivement par $\pgcd(p(x),p'(x))$, on peut supposer que le polynôme $p(x)$ ne possède que des facteurs simples. On obtiendra alors pour $A$ une matrice flèche-diagonale (diagonale par blocs dont les blocs sont des matrices flèches).
\air
On suppose que le polynôme s'écrit toujours 
$$p(x)=\prod_{i=1}^r(x-\alpha_i)\prod_{j=1}^s(-(x-\beta_j)^2-\gamma_j^2)$$ 
même si on n'a plus accès aux réels $\alpha_i,\beta_j$ et $\gamma_j\not =0$.

\air{\bf Localisation des racines}\par
Par exemple grâce aux suites de Sturm, on peut déterminer le nombre $r$ de racines réelles de $p(x)$ et aussi les localiser.
On supposera $\alpha_1<\ldots<\alpha_r$. 

\air{\bf Entrelacement}\par
On choisit une famille arbitraire de réels $\lambda_1,\ldots,\lambda_{r-1}$ qui entrelacent les racines réelles, c'est-à-dire telle que
$$\alpha_1<\lambda_1<\alpha_2<\lambda_2<\ldots<\alpha_{r-1}<\lambda_{r-1}<\alpha_r$$

On cherche alors la matrice $A$ sous la forme flèche suivante :
$$A=
\left(\begin{array}{ccccccccc}
\mu_1&\nu_1&&&&&&&k_1\\
\nu_1&-\mu_1&&&&&&&l_1\\
&&\ddots&&&&&&\vdots\\
&&&\mu_s&\nu_s&&&&k_s\\
&&&\nu_s&-\mu_s&&&&l_s\\
&&&&&\lambda_1&&&h_1\\
&&&&&&\ddots&&\vdots\\
&&&&&&&\lambda_{r-1}&h_{r-1}\\
k_1&l_1&\ldots&k_s&l_s&h_1&\ldots&h_{r-1}&e\\
\end{array}\right)
$$

où $(\mu_j,\nu_j)\in\RR\times\RR^*$ pour $j=1\ldots s$ sont des réels arbitraires et les $h_i,k_j,l_j$ et $e$ sont des réels à déterminer pour $i=1\ldots r$ et $j=1\ldots s$.

\air{\bf Interpolation}\par 
En se servant \ref{schur} (ii), on calcule 
$$\begin{array}{rcl}
q(x)&=&\det(xJ-A)\\
&=&\left(\begin{array}{ccccccccc}
x-\mu_1&-\nu_1&&&&&&&-k_1\\
-\nu_1&-x+\mu_1&&&&&&&-l_1\\
&&\ddots&&&&&&\vdots\\
&&&x-\mu_s&-\nu_s&&&&-k_s\\
&&&-\nu_s&-x+\mu_s&&&&-l_s\\
&&&&&x-\lambda_1&&&-h_1\\
&&&&&&\ddots&&\vdots\\
&&&&&&&x-\lambda_{r-1}&-h_{r-1}\\
-k_1&-l_1&\ldots&-k_s&-l_s&-h_1&\ldots&-h_{r-1}&x-e\\
\end{array}\right)\\
&&\\
&=&(x-e)u(x)v(x)-v(x)\sum_{i=1}^{r-1}h_i^2u_i(x)-u(x)\sum_{j=1}^s((-x+\mu_j)k_j^2+2\nu_jk_jl_j+(x-\mu_j)l_j^2)v_j(x)
\end{array}$$
où\par 
$\left\{\begin{array}{ccl}
u(x)&=&\prod_{i=1}^{r-1}(x-\lambda_i)\\
u_i(x)&=&\frac{u(x)}{x-\lambda_i}\\
v(x)&=&\prod_{j=1}^s-((x-\mu_j)^2+\nu_j^2)\\
v_j(x)&=&\frac{v(x)}{-((x-\mu_j)^2+\nu_j^2)}
\end{array}\right.$\air

Du coup, $q(\lambda_i)=-v(\lambda_i)h_i^2u_i(\lambda_i)$. Or $p(\lambda_i)=\prod_{j=1}^s-((\lambda_i-\beta_j)^2+\gamma_j^2)\times\prod_{k=1}^r(\lambda_i-\alpha_k)$, dont le premier produit est du signe de $v(\lambda_i)$ (qui est aussi celui de $(-1)^s$). Ainsi, grâce à la condition d'entrelacement, il est possible de trouver un réel $h_i$ tel que $q(\lambda_i)=p(\lambda_i)$ pour tout $i=1\ldots r-1$.\par
Puis 
$$
\begin{array}{ccl}
q(\mu_j+I\nu_j)&=&-u(\mu_j+I\nu_j)((-I\nu_j)k_j^2+2\nu_jk_jl_j+(I\nu_j)l_j^2)v_j(\mu_j+I\nu_j)\\
&=&-u(\mu_j+I\nu_j)(-I\nu_j)(k_j+Il_j)^ 2 v_j(\mu_j+I\nu_j)
\end{array}$$

Il est donc possible de trouver $(k_j,l_j)\in\RR^2$ tels que $q(\mu_j+I\nu_j)=p(\mu_j+I\nu_j)$ pour tout $j=1\ldots s$.\par
Pour finir, si on choisi $e$ tel que $(-1)^s(\sum_{j=1}^s2\mu_j-\sum_{i=1}^{r-1}{\lambda_i}-e)$ vaut le coefficient en $x^{d-1}$ de $p(x)$, on est assuré que $p(x)=q(x)$.

\air{\bf Cas particulier}\par 
Nous nous devons traiter à part le cas où $p(x)$ ne possède {\it aucune} racine réelle. Dans ce cas, 
on pose 
$$J=\left(\begin{array}{ccccccc}
1&0&&&&&\\
0&-1&&&&&\\
&&\ddots&&&&\\
&&&1&0&&\\
&&&0&-1&&\\
&&&&&-1&\\
&&&&&&1\\
\end{array}\right)=\left(\Id_1\bigoplus(-\Id_1)\right)^{s-1}\bigoplus\left((-\Id_1)\bigoplus \Id_1\right)$$ 
et on cherche $A$ sous la forme 
$$A=
\left(\begin{array}{ccccccc}
\mu_1&\nu_1&&&&&k_1\\
\nu_1&-\mu_1&&&&&l_1\\
&&\ddots&&&\\
&&&\mu_{s-1}&\nu_{s-1}&&k_{s-1}\\
&&&\nu_{s-1}&-\mu_{s-1}&&l_{s-1}\\
&&&&&\lambda&h\\
k_1&l_1&\ldots&k_{s-1}&l_{s-1}&h&e
\end{array}\right)
$$
En fait, le seul aménagement du cas précédent consiste à remarquer que la quantité $$(\det(xJ-A))_{(x=-\lambda)}=-h^2\times\prod_{j=1}^{s-1}(-(-\lambda-\mu_j)^2-\nu_j^2)$$ est du même signe que $p(-\lambda)$, ce qui permet de trouver un réel $h$ qui assure l'interpolation au point $x=-\lambda$. Le reste de l'argument est analogue.    
\end{preuve}


\section{Racines réelles des représentations déterminantales}
Intéressons-nous dorénavant à la propriété réciproque du Théorème \ref{effective_rep_determ_fleche}.

\begin{thm}\label{recip}
Soit $A$ une matrice symétrique de $\SS\RR^{d\times d}$ avec $d=r+2s$. Alors, le polynôme $p(x)=\det\left(xJ_s-A\right)$, où $J_s=\left(\Id_{r+s}\bigoplus(-\Id_s)\right)$, possède au moins $r$ racines réelles comptées avec multiplicité.
\end{thm}
\begin{preuve}
Si on reformule le problème en terme d'endomorphisme auto-adjoint relativement à la forme bilinéaire symétrique donnée par la matrice $J_s$, on est ramené à montrer le Théorème \ref{spectre_autoadj}.
En effet, si on considère que $J_s$ est la matrice d'une forme bilinéaire symétrique (de signature $(r+s,s)$) sur $\RR^{r+2s}$, alors relativement à cette forme bilinéaire symétrique, l'adjointe $B^{(*_s)}$ de la matrice $B$ est donnée par la formule $B^{(*_s)}=J_sB^TJ_s$. 
Il reste alors à remarquer que $p(x)=(-1)^s\det(x\Id-AJ_s)$ et que $AJ_s$ est $J_s$-autoadjointe.
\end{preuve}

Avant d'énoncer le résultat proprement dit, notons que si on se donne une forme quadratique de signature $(r+s,s)$ sur $\RR^{r+2s}$ par une matrice symétrique $S$, alors la forme bilinéaire associée s'écrit $\langle X,Y\rangle=X^TSY$. Or on peut décomposer $S=Q^TJ_sQ$ avec $Q$ inversible, ainsi $\langle X,Y\rangle=(QX)^TJ_s(QY)$, donc à un changement de coordonnées près, on peut supposer que la forme quadratique est donnée par la matrice $J_s$. 

\begin{thm}\label{spectre_autoadj}
Dans l'espace $\RR^n$, muni d'une forme quadratique de signature 
$(r+s,s)$, tout endomorphisme auto-adjoint possède au moins $r$ racines réelles comptées avec multiplicités.
\end{thm} 

\begin{preuve}
On suppose donné l'endomorphisme auto-adjoint par sa matrice $A$ dans la base canonique. La démonstration se fait par récurrence sur $s$. On commencera par traiter les cas $s=0$ et $s=1$ pour fixer les idées. Ensuite, on s'attaquera à l'hérédité dans le cas général. L'idée consiste grosso-modo en une orthonormalisation relativement à une forme bilinéaire symétrique de type $J_k$ qui fera apparaître une matrice flèche dont on saura calculer le déterminant. On pourra alors en déduire le nombre de racines souhaité. En fait, l'étape de $J_k$-orthonormalisation ne marche bien que si la matrice considérée n'a que des valeurs propres simples. Pour se ramener à cette situation, on ajoutera de nouvelles variables infinitésimales.\air

\noindent{\bf Quelques notations}\par 
Soient $(r,s)\in\N^2$ et pour tout entier $k=0\ldots s$, notons $n_k=r+s+k$ et $J_k=\Id_{r+s}\bigoplus(-\Id_k)\in\RR^{n_k\times n_k}$.

\par
La matrice $A$ est de la forme
$$A=\left(\begin{array}{cc}
A_0&-H\\
H^T&D
\end{array}\right)$$
avec $H\in\RR^{n_0\times s}$, $D=(d_{i,j})\in\SS\RR^{s\times s}$ et $A_0\in\SS\RR^{n_0\times n_0}$.\par 

Appelons $A_k$ la $n_k$-ième matrice extraite principale de $A$,
de sorte que, pour tout $k=0\ldots s-1$, on a des matrices $H_{k}\in\RR^{n_k\times 1}$ et $L_{k}\in\RR^{1\times n_k}$ telles que  :
$$A_{k+1}=\left(\begin{array}{cc}
A_k&-H_{k+1}\\
L_{k+1}&d_{k+1,k+1}\\
\end{array}\right)$$

On introduit le polynôme caractéristique $p(x)=p_A(x)=\det(x\Id-A)$ de $A$, et il s'agit de montrer qu'il possède au moins $r$ valeurs propres réelles, comptées avec multiplicité. Nous introduisons aussi les polynômes caractéristiques $p_k(x)=\det(x\Id_{n_k}-A_k)$. 

\par

\air On procède par récurrence sur $s$. \air
\noindent{\bf Les cas $s=0$ et $s=1$}\par 
Si $s=0$, le résultat est clair : c'est le théorème spectral usuel pour les matrices symétriques réelles. En vue de préparer l'étape suivante, explicitons son utilisation.\par  Il existe une matrice orthogonale $Q_0\in\RR^{n_0\times n_0}$, c'est-à-dire telle que $Q_0^TQ_O=\Id$ (on dira aussi qu'elle est $J_0$-orthonormale, autrement dit que $Q_0^{(*_0)}Q_0=\Id$ où l'adjoint ${}^{(*_0)}$ est relatif à la forme bilinéaire symétrique donnée par $J_0$) telle que la matrice $B_0=Q_0^TA_0Q_0$ est diagonale : 

$$B_0=\left(\begin{array}{ccc}
\lambda_1&0&0\\
0&\ddots&0\\
0&0&\lambda_{n_0}
\end{array}\right).$$

Posons  
$$A'_{0}=B_0+\left(\begin{array}{cccc}
\epsilon_0&0&0&0\\
0&2\epsilon_0&0&0\\
0&0&\ddots&0\\
0&0&0&n_0\epsilon_0
\end{array}\right).$$
On obtient alors une matrice $A'_{0}\in\left(\RR\langle\epsilon_0\rangle_0\right)^{n_0\times n_0}$, diagonale, dont toutes les valeurs propres sont simples, et telle que $\lim_{\epsilon_0\to 0}A'_{0}=B_0=Q_0^TA_0Q_0$.
\air

Traitons aussi le cas $s=1$, pour fixer les idées, même si l'hérédité dans le cas général comprendra ce cas là. 
\par

Du fait que $A_1=\left(\begin{array}{cc}
A_0&-H_1\\
L_1&d_{1,1}
\end{array}\right)$ avec $H_1^T=L_1$, on introduit la matrice $B_1\in\R^{n_1\times n_1}$ telle que $$B_1=\left(Q_0\bigoplus\Id\right)^{(*_1)} A_1\left(Q_0\bigoplus\Id\right) =
\left(\begin{array}{cc}
B_0&-H'_1\\
L'_1&d'_{1,1}
\end{array}\right)$$ où ${}^{(*_1)}$ est l'adjoint relativement à $J_1$. Notons au passage que $\left(Q_0\bigoplus\Id\right)^{(*_1)}=Q^{(*_0)}_0\bigoplus\Id_1$.\air 

On peut écrire, par le Lemme \ref{schur} (i) :  
$$p_1(x)=\det(x\Id_{n_1}-A_1)
=\det\left(
\begin{array}{cc}
x\Id_{n_0}-B_0&H'_1\\
-L'_1&x-d'_{1,1}
\end{array}
\right)=
(x-d'_{1,1})p_0(x)+\sum_{i=1}^{n_0}h_i^2\prod_{j\not =i}(x-\lambda_j)
$$\par
avec $(H'_1)^T=(h_1,\ldots,h_{n_0})=L'_1$.\air
Supposons, dans un premier temps, que 
\ib
\item[(i)] les $\lambda_i$ sont distincts : mettons $\lambda_1<\ldots<\lambda_{n_0}$,
\item[(ii)] les $h_i$ sont non nuls. 
\ie
En évaluant en $\lambda_i$, et en examinant le signe de $p_1(\lambda_i)$ pour $i=1\ldots n_0$, on remarque que le polynôme $p_1(x)$ possède $n_0-1$ racines réelles distinctes (entrelacées avec les $\lambda_i$), ce qui donne le résultat souhaité.\air

Pour traiter le cas où les conditions $(i)$ ou $(ii)$ sont mises en défaut, il suffit de changer de corps réel clos de base. On fait appel alors aux techniques employées dans le Lemme \ref{approx_schur}.\par

On est alors conduit à introduire la matrice de $\R\langle\langle\epsilon_0,\epsilon_1,\theta_1\rangle\rangle^{n_1\times n_1}$ suivante : 
$$A'_1=
\left(\begin{array}{cc}
A'_0&-H'_1\\
L'_1&d'_1
\end{array}\right)+\left(\begin{array}{ccccc}
0&\ldots&0&-\epsilon_1\\
\vdots&&\vdots&\vdots\\
0&\ldots&0&-\epsilon_1\\
\epsilon_1&\ldots&\epsilon_1&\theta_1
\end{array}\right)$$


On obtient alors une nouvelle matrice $A_1'$ qui possède, d'après le traitement ci-dessus du cas générique, au moins $n_0-1$ racines (distinctes) dans $R\langle\langle\epsilon_0,\epsilon_1,\theta_1\rangle\rangle$. Comme Ainsi $\lim_{(\epsilon_0,\epsilon_1,\theta_1)\to 0}A_1'=B_1$ est semblable à $A_1$, alors $A_1$ possède au moins $n_0-1$ racines réelles, comptées avec multiplicité.\par 

Notons que l'ajout de la variable $\theta_1$ n'a pas servi jusqu'alors, elle permet cependant  de préparer l'étape suivante de la récurrence.
Elle sert à affirmer que, d'après le Lemme \ref{approx_schur}, toutes les valeurs propres de $A'_1$ sont distinctes dans $R\langle\langle\epsilon_0,\epsilon_1,\theta_1\rangle\rangle[I]$. De plus, on peut garder à l'esprit que $B_1=Q_1^{(*_1)}A_1Q_1\in \R^{n_1\times n_1}$, avec $Q_1=\left(Q_0\bigoplus\Id_1\right)$.\par 


\air
{\bf Hérédité dans le cas général\par} 
Supposons donnée une matrice $A'_k\in \SS R_k^{n_k\times n_k}$ telle que
\ib
\item[(i)] On dispose d'une suite d'extensions réelles closes $\RR\rightarrow R_1\rightarrow\ldots\rightarrow R_k$ avec $R_{i+1}=R_i\langle\langle\epsilon_{i+1},\theta_{i+1}\rangle\rangle$ pour tout $i=1\ldots k-1$, 
\item[(ii)] $A'_k$ ne possède que des valeurs propres simples dans $R_k[I]$, 
\item[(iii)] $A'_k$ est $J_k$-autoadjointe, c'est-à-dire $A'_k=J_k(A'_k)^TJ_k$,
\item[(iv)] $A'_k$ est à coefficients dans $\R\langle\epsilon_0,\epsilon_1,\theta_1,\ldots,\epsilon_k,\theta_k\rangle_0$ et on a $\lim_{(\epsilon_0,\epsilon_1,\theta_1,\ldots,\epsilon_k,\theta_k)\to 0}A'_k=Q_k^{(*_k)}A_kQ_k$, avec $Q_k$ matrice $J_k$-orthonormale de $\RR^{n_k\times n_k}$ (i.e.  $Q_k^{(*_k)}Q_k=\Id$),
\item[(v)] $A'_k$ possède au moins $n_0-k$ valeurs propres (distinctes) dans $R_k$.
\ie
\par
Avant de commencer proprement dit la démonstration, établissons un Lemme préliminaire d'algèbre bilinéaire.

\begin{lem}\label{autoadj}
Notons $\langle\cdot,\cdot\rangle$ la forme bilinéaire symétrique sur $R^n$ donnée par la matrice $J=\Id_p\bigoplus(-\Id_{q})$, où $p+q=n$ et $R$ est réel clos. On considère $A\in R^{n\times n}$ une matrice $J$-autoadjointe, c'est-à-dire pour tous $X,Y\in R^{n\times 1}$, $\langle AX,Y\rangle=\langle X,AY\rangle$. On note $A^*$ la matrice adjointe de $A$ relativement à $J$.
\ib
\item[1)] La matrice $A$ est $J$-autoadjointe si et seulement si l'une des conditions équivalentes suivantes est satisfaite
\ib
\item[(i)] $AJ$ est symétrique, 
\item[(ii )] $A=BJ$ avec $B$ symétrique
\item[(iii)] $A=\left(\begin{array}{cc}
U&-V\\
V^T&W
\end{array}\right)$ avec $U\in \SS R^{p\times p}$, $V\in R^{p\times q}$, $W\in \SS R^{q\times q}$.
\ie
\item[2)] Si $u$ et $v$ sont deux vecteurs propres de $A$ associés respectivement à deux valeurs propres distinctes $\lambda$ et $\mu$ dans $R[I]$, alors $u$ et $v$ sont $J$-orthogonaux.
\item[3)] Si $A$ possède $n$ valeurs propres distinctes dans $R[I]$, alors il existe $P\in R^{n\times n}$ une matrice $J$-orthonormale (c'est-à-dire telle que $P^*P=\Id_n$) telle que $P^*AP$ soit, à permutation des coordonnées près, une matrice de $R^{n\times n}$ diagonale par blocs, avec des blocs de taille $1$ ou de taille $2$ sans valeur propre dans $R$. 
\ie 
\end{lem}
\begin{preuve}
1) Pour les points $(i)$ et $(ii)$, il suffit de noter que $A^T=JAJ$ équivaut à $(AJ)^T=AJ$. La dernière assertion résulte juste d'un calcul matriciel par bloc : \par si $A=\left(\begin{array}{cc}U&-V\\V'&W\end{array}\right)$ alors $JAJ=\left(\begin{array}{cc}U&V\\-V'&W\end{array}\right)$.

2) Si $\lambda$ et $\mu$ sont dans $R$, il suffit d'écrire $\langle Au,v\rangle=\lambda\langle u,v\rangle=\langle u,Av\rangle=\mu\langle u,v\rangle$. Comme $\lambda\not =\mu$, on a forcément $\langle u,v\rangle=0$.\par
Si $\lambda$ ou $\mu$ sont dans $R[I]$, on étend la forme bilinéaire symétrique $\langle\cdot,\cdot\rangle$ à $R[I]^{n}$, ce qui permet d'écrire les mêmes égalités que précédemment. 

3) {\`A} une valeur propre $\lambda\in R$, on associe un vecteur propre $u$ et à une valeur propre $\mu+I\nu\in R[I]\setminus R$, on associe un vecteur propre $v+Iw$, avec $v,w\in R^n$. Du coup $\mu-I\nu$ est associée $v-Iw$.\par
On dispose alors d'une base de $R^n$ formée des vecteurs $(u_1,\ldots,u_r,v_1,w_1,\ldots,v_s,w_s)$. En vertu de $1)$, on a pour tout $i=1\ldots r$, $u_i^\perp\supset\vect(u_1,\ldots,u_{i-1},u_{i+1},\ldots, u_r,v_1,w_1,\ldots,v_s,w_s)$, l'orthogonalité considérée est bien sûr la $J$-orthogonalité dans $R^n$ relativement à la forme bilinéaire symétrique $\langle\cdot,\cdot\rangle$. Par considération des dimensions, l'inclusion précédente est une égalité.\par
De même, pour $j=1\ldots s$, $v_j^\perp$ et $w_j^\perp$ contiennent tout deux l'espace vectoriel $$\vect(u_1,\ldots,u_r,v_1,w_1,\ldots,v_{j-1},w_{j-1},v_{j+1},w_{j+1}\ldots,v_s,w_s).$$ Et on a aussi $\langle v_j,v_j\rangle+\langle w_j,w_j\rangle=0$ par $J$-orthogonalité de $v_j+Iw_j$ et de $v_j-Iw_j$.\air
Trois situations se présentent :
\ib
\item[a)] si $\langle v_j,w_j\rangle=0$, alors $v_j^\perp$ contient aussi $w_j$ et $w_j^\perp$ contient aussi $v_j$. En particulier, par considération des dimensions de $v_j^\perp $ et $w_j^\perp $, on déduit que $v_j$ et $w_j$ ne sont pas isotropes.
\item[b)] si $\langle v_j,w_j\rangle\not =0$ et $v_j$ est isotrope, alors $w_j$ est aussi isotrope et on pose par exemple $v'_j=v_j+w_j$ et $w'_j=v_j-w_j$. Ainsi $\langle v'_j,w'_j\rangle=0$ et en remplaçant la famille $(v_j,w_j)$ par la famille $(v'_j,w'_j)$ on obtient les mêmes conclusions que précédemment.
\item[c)] si $\langle v_j,w_j\rangle\not =0$ et $v_j$ n'est pas isotrope, alors on remplace $(v_j,w_j)$ par $(v_j,v_j-\frac{\langle v_j,v_j\rangle}{\langle v_j,w_j\rangle}w_j)$ et on  obtient les mêmes conclusions que précédemment.
\ie
Dans tous les cas, on obtient une base $J$-orthogonale, chacun des vecteurs de la base étant nécessairement non isotrope. On la normalise, de sorte que pour chaque vecteur $u$ de la nouvelle base ${\cal B}$, on ait $\langle u,u\rangle=\pm 1$. D'après le Théorème d'inertie de Sylvester, on compte exactement $p$ fois $+1$ et $n-p$ fois $-1$. Quitte à permuter les vecteurs, on peut supposer que la matrice de Gram associée à la base ${\cal B}$ est la matrice $J$. Notons $P$ la matrice de passage de la base canonique à la base ${\cal B}$. \par Nous avons bien $P^*P=\Id$ et $D=PAP^*$ qui est, à permutation près des coordonnées, une matrice diagonale par blocs dont chaque bloc est de taille $1$ ou $2$.\par

Un bloc $M_j$ de taille $2$ correspond à la matrice de la restriction à un plan vectoriel du type $\vect(v_j,w_j)$ de l'endomorphisme canoniquement associé à $A$. Ainsi ce bloc $M_j$ a deux valeurs propres distinctes dans $R[I]$, à savoir $\mu_j\pm\nu_j$. 
\end{preuve}

\begin{rem}
D'après le Lemme \ref{autoadj} 1), et comme la matrice $D$ est encore $J$-autoajointe, elle est de la forme  $D=\left(\begin{array}{cc}
U&-V\\
V^T&W
\end{array}\right)$ avec $U\in \SS R^{p\times p}$, $V\in R^{p\times q}$, $W\in \SS R^{q\times q}$.
Donc pour tout $j$, les vecteurs $v_j$ et $w_j$ ne peuvent pas être simultanément en une position indexée par $\{1\ldots p\}$, ni être simultanément en une position indexée par $\{p+1\ldots p+q\}$ dans la nouvelle base $J$-orthonormale {\cal B}, sinon cela signifierait que $M_j$ est symétrique, une contradiction.
\end{rem}
\air

Venons-en au coeur de la démonstration de l'hérédité :\par
Du fait que $A_{k+1}=\left(\begin{array}{cc}
A_k&-H_{k+1}\\
L_{k+1}&d_{k+1}
\end{array}\right)$, 
on pose $$B_{k+1}=\left(Q_k\bigoplus\Id_1\right)^{(*_{k+1})} A_k\left(Q_k\bigoplus\Id_1\right) =
\left(\begin{array}{cc}
Q_k^{(*_k)}A_kQ_k&-H'_{k+1}\\
L'_{k+1}&d'_{k+1}
\end{array}\right)\in \RR^{n_k\times n_k}$$ 

Notons qu'on a $\left(Q_k\bigoplus\Id_1\right)^{(*_{k+1})}=Q_k^{(*_k)}\bigoplus\Id_1$. \par

Introduisons ensuite la matrice de $R_{k+1}^{n_{k+1}\times n_{k+1}}$ suivante : 
$$B'_{k+1}=
\left(\begin{array}{cc}
A'_k&-H'_{k+1}\\
L'_{k+1}&d'_{k+1}
\end{array}\right)$$


D'après le Lemme \ref{autoadj} et les points $(ii)$ et $(iii)$ de l'hypothèse de récurrence, il existe $P_k\in R_k^{n_k\times n_k}$ une matrice $J_k$-orthonormale, telle que $P_k^*A'_kP_k$ soit, à permutation des coordonnées près, diagonale par blocs, avec $a$ blocs de taille $1$ et $b$ blocs $M_1,\ldots, M_{b}$ de taille $2$, dont les polynômes caractéristiques sont toujours strictement positifs. On a les inégalités $a\geq n-k$ et $b\leq k$.\air
Posons alors $$\widetilde{A_{k+1}}=\left(P_k\bigoplus\Id_1\right)^{(*_{k+1})}B'_{k+1} \left(P_k\bigoplus\Id_1\right)$$
et finalement
$$A'_{k+1}=\widetilde{A_{k+1}}+\epsilon_{k+1}U_{k+1}+\theta_{k+1}V_{k+1}$$ avec $V_{k+1}=\left(\begin{array}{cccc} 0&\ldots&0&0\\ \vdots&&\vdots&\vdots\\ 0&\ldots&0&0\\ 0&\ldots&0&1 \end{array}\right)$ et $U_{k+1}=\left(\begin{array}{cccccc} 
0&\ldots&\ldots&\ldots&0&-1\\ 
\vdots&&&&\vdots&\vdots\\ 
\vdots&&&&\vdots&-1\\ 
\vdots&&&&\vdots&0\\ 
0&\ldots&\ldots&\ldots&0&\vdots\\ 
+1&\ldots&+1&0&\ldots&0 \end{array}\right)$, cette dernière matrice comptant exactement $k+1$ coefficients nuls sur la dernière colonne et la dernière ligne.\air

Comme $P_k$ est $J_k$-orthonormale, alors $\left(P_k\bigoplus\Id_1\right)$ est $J_{k+1}$-orthonormale, et du fait que $A_{k+1}$ est $J_{k+1}$-autoadjointe, il en sera de même de $\widetilde{A_{k+1}}$ et donc de $A'_{k+1}$ qui est donc de la forme 

$$A'_{k+1}=\left(\begin{array}{cc}
P_k^{(*_k)}A'_kP_k&-H''_{k+1}\\
L''_{k+1}&d''_{k+1,k+1}
\end{array}\right)$$\air

Par construction, $$\lim_{(\epsilon_0,\epsilon_1,\theta_1,\ldots,\epsilon_{k+1},\theta_{k+1})\to 0} A'_{k+1}=Q^{(*_{k+1})}_{k+1}A_{k+1}Q_{k+1}$$ avec $Q_{k+1}=(P_k\bigoplus \Id_1)\times(Q_k\bigoplus \Id_1)\in\R^{n_{k+1}\times n_{k+1}}$.
\air
On a donc construit une extension réelle close $R_k\rightarrow R_{k+1}=R_k\langle\langle\epsilon_{k+1},\theta_{k+1}\rangle\rangle$ et une matrice $A'_{k+1}\in R_{k+1}^{n_{k+1},n_{k+1}}$ satisfaisant les points $(i),(ii)$ (grâce au Lemme \ref{approx_schur}),$(iii),(iv)$. Le point $(v)$ découlera de l'étude du polynôme $p_{A'_{k+1}}(x)$ que l'on conduit en examinant les deux cas suivants (on aurait pu se limiter au second cas mais le traitement du premier cas éclaire la démarche).

\ib
\item[$\bullet$]
Supposons que $k=b$.
D'après le Lemme \ref{autoadj} 1) et aussi d'après la remarque suivant sa démonstration, on peut supposer, à une permutation des coordonnées près, que :

$$A'_{k+1}=\left(\begin{array}{ccccccc}
\lambda_1&&&&&&-h_1-\epsilon_{k+1}\\
&\ddots&&&&&\vdots\\
&&\lambda_a&&&&-h_a-\epsilon_{k+1}\\
&&&M_1&&&\left(\begin{array}{c}-g_1-\epsilon_{k+1}\\ l_1\end{array}\right)\\
&&&&\ddots&&\vdots\\
&&&&&M_b&\left(\begin{array}{c}-g_b-\epsilon_{k+1}\\ l_b\end{array}\right)\\
h_1+\epsilon_{k+1}&\ldots&h_a+\epsilon_{k+1}&(g_1+\epsilon_{k+1},l_1)&\ldots&(g_b+\epsilon_{k+1},l_b)&d+\theta_{k+1}
\end{array}\right)$$
où $h_i,g_j,l_j\in R_k$ pour $i=1\ldots a$, $j=1\ldots b$.\air

On fait alors appel au Lemme \ref{schur}(ii), et par un calcul analogue à celui de la démonstration de \ref{effective_rep_determ_fleche}, on obtient
$$p_{A'_{k+1}}(x)=(x-d-\theta_{k+1})p_{A'_k}(x)+\left(\sum_{i=1}^{r_k}(h_i+\epsilon_{k+1})^2u_i(x)\right)v(x)+u(x)w(x)$$
avec $u(x)=\prod_{i=1}^{a}(x-\lambda_i)$, $u_i(x)=\frac{u(x)}{x-\lambda_i}$, $v(x)=\prod_{j=1}^b p_{M_j}(x)$ et $w(x)$ est un polynôme qu'on n'a pas besoin d'expliciter ici.\par
En évaluant successivement en $\lambda_i$ pour $i=1\ldots a$, et en utilisant le fait que $v(x)$ est toujours strictement positif et que $h_i+\epsilon_{k+1}\not=0$, on obtient une suite d'éléments $p_{A'_{k+1}}(\lambda_1),\ldots,p_{A'_{k+1}}(\lambda_{r_k})$ qui présente $a$ changements de signes, donc $p_{A_{k+1}}(x)$ possède au moins $a-1=n_0-(k+1)$ racines (distinctes) dans $R_{k+1}$.
\air

\item[$\bullet$] Supposons maintenant que $k>b$. On pose $c=n_0-b$ et $c+d=a$ (du coup $d=k-b$).\par

{\`A} une permutation près des coordonnées, on a l'égalité\vskip0,5cm 
{\footnotesize $$A'_{k+1}=\left(\begin{array}{ccccccccccc}
\lambda_1&&&&&&&&&-h_1-\epsilon_{k+1}\\
&\ddots&&&&&&&&\vdots\\
&&\lambda_c&&&&&&&-h_c-\epsilon_{k+1}\\
&&&\mu_{1}&&&&&&m_{1}\\
&&&&\ddots&&&&&\vdots\\
&&&&&\mu_d&&&&m_d\\
&&&&&&M_1&&&\left(\begin{array}{c}-g_1-\epsilon_{k+1}\\ l_1\end{array}\right)\\
&&&&&&&\ddots&&\vdots\\
&&&&&&&&M_b&\left(\begin{array}{c}-g_b-\epsilon_{k+1}\\ l_b\end{array}\right)\\
h_1+\epsilon_{k+1}&\ldots&h_c+\epsilon_{k+1}&m_{1}&\ldots&m_d&(g_1+\epsilon_{k+1},l_1)&\ldots&(g_b+\epsilon_{k+1},l_b)&d+\theta_{k+1}
\end{array}\right)$$}\vskip0,5cm
où $h_i,g_j,l_j,m_v\in R_k$ pour $i=1\ldots c$, $j=1\ldots b$, $t=1\ldots d$. On suppose aussi par commodité que $\lambda_1<\ldots<\lambda_c$.\par 
Etudions le nombre de changements de signes de la suite $p_{A'_{k+1}}(\lambda_1),\ldots,p_{A'_{k+1}}(\lambda_c)$ sachant que $p_{A'_{k+1}}(\lambda_i)$ est du signe de $$\prod_{j=1,j\not = i}^c(\lambda_i-\lambda_j)\prod_{t=1}^d(\lambda_i-\mu_t).$$ 
Si on pose $u(x)=\prod_{i=1}^c(x-\lambda_i)$ et $u_i(x)=\frac{u(x)}{x-\lambda_i}$, alors la suite $u_1(\lambda_1),\ldots,u_c(\lambda_c)$ présente $c-1$ changement de signes, puis la suite $u_1(\lambda_1)(\lambda_1-\mu_1),\ldots,u_c(\lambda_c)(\lambda_c-\mu_1)$ présente au moins $c-2$ changements de signes, et ainsi de suite jusqu'à la suite $u_1(\lambda_1)\prod_{t=1}^d(\lambda_1-\mu_t),\ldots,u_c(\lambda_c)\prod_{t=1}^d(\lambda_c-\mu_t)$ qui présente au moins $c-d-1$ changements de signes. Au final, $p_{A'_{k+1}}(x)$ possède au moins $n_0-k-1$ racines distinctes dans $R_{k+1}$.
\ie

\par
Ce qui achève la démonstration de l'hérédité.\air

Au bout du compte, on parvient à une matrice $A'_s$ qui possède au moins $n_s$ valeurs propres simples dans $\RR\langle\langle\epsilon_0,\epsilon_1,\theta_1,\ldots,\epsilon_s\theta_s\rangle\rangle$ et telle que $\lim_{(\epsilon_0,\epsilon_1,\theta_1,\ldots,\epsilon_s\theta_s)\to 0}A'_s=A_s$. Donc $A_s=A$ a au moins $n_0-s=r$ valeurs propres réelles comptées avec multiplicité.
\end{preuve}

\begin{rem}
On peut rapprocher le principe de cette démonstration de la démonstration du théorème spectral classique pour les matrices symétriques réelles donnée dans \cite{BPR}.  
\end{rem}



\section{D'une représentation flèche-diagonale à une autre}
On peut affiner le résultat \ref{effective_rep_determ_fleche}, qui conjointement avec \ref{recip} peut se reformuler de la sorte :\par

\begin{thm}
Tout polynôme $p(x)\in\RR[x]$ de degré $d=r+2s$, possède au moins $r$ racines réelles comptées avec multiplicité si et seulement s'il admet une représentation effective de signature $(r+s,s)$.
\end{thm}

\begin{preuve}
Il suffit de montrer que si le polynôme $p(x)$ admet une représentation déterminantale effective de signature $(r+s,s)$ avec $r\geq 1$, alors il en admet aussi une de signature $(r+s-1,s+1)$.\air

Après application de  \ref{recip}, puis de \ref{effective_rep_determ_fleche}, cette propriété (sans tenir compte du caractère d'effectivité), se ramène à la propriété évidente : si le polynôme $p(x)$ possède au moins $r$ racines, alors il possède au moins $r-1$ racines.\par
Voyons maintenant ce que l'ont peut dire du point de vue effectif. En général, il n'est pas évident de passer formellement d'une représentation déterminantale $p(x)=\det(J-xA)$ de signature $(r+s,s)$ à une représentation déterminantale $p(x)=\det(J'-xA')$ de signature $(r+s-1,s+1)$. Par contre, dans le cas d'une représentation flèche-diagonale, le procédé est aisé et se résume encore une fois au calcul de déterminant donné dans le Lemme \ref{schur}.\par
De nouveau, on supposera que $p(x)$ n'a que des facteurs simples et on préférera travailler avec le polynôme $p^*(x)=\det(xJ-A)$.\air  
\ib
\item[a)] Si $r=1$, il suffit de multiplier par $-\Id$, c'est-à-dire $p(x)=(-1)^d\det(x(-J)-(-A))$ et $-J$ est bien de signature $(s,s+1)$.
\item[b)] Supposons que $r=2$ et que les matrices $A$ et $J$ sont 
données par 
$$A=
\left(\begin{array}{ccccccccc}
\mu_1&\nu_1&&&&&k_1\\
\nu_1&-\mu_1&&&&&l_1\\
&&\ddots&&&&\vdots\\
&&&\mu_s&\nu_s&&k_s\\
&&&\nu_s&-\mu_s&&l_s\\
&&&&&\lambda&h\\
k_1&l_1&\ldots&k_s&l_s&h&e\\
\end{array}\right)
$$
et 

$$J=
\left(\begin{array}{ccccccc}
1&&&&&&\\
&-1&&&&&\\
&&\ddots&&&&\\
&&&1&&&\\
&&&&-1&&\\
&&&&&1&\\
&&&&&&1\\
\end{array}\right)=\left(\Id_1\bigoplus(-\Id_1)\right)^s\bigoplus\Id_2
$$

Posons alors $J'=(\Id_1\bigoplus(-\Id_1))^{s}\bigoplus(-\Id_{1})\bigoplus\Id_1$ et 

$$A'=
\left(\begin{array}{ccccccc}
\mu_1&\nu_1&&&&&k'_1\\
\nu_1&-\mu_1&&&&&l'_1\\
&&\ddots&&&&\vdots\\
&&&\mu_{s}&\nu_{s}&&k'_{s}\\
&&&\nu_{s}&-\mu_{s}&&l'_{s}\\
&&&&&-\lambda&h'\\
k'_1&l'_1&\ldots&k'_{s}&l'_{s}&h'&e'\\
\end{array}\right)
$$

D'après les formules de la démonstration du Théorème \ref{effective_rep_determ_fleche}, il suffit de prendre les réels $h',k'_j,l'_j$ pour $j=1\ldots s$ et $e'$ tels que :
$$\left\{\begin{array}{ccc}
h'&=&h\\
k'_j&=&l_j\\
l'_j&=&-k_j\\
e'&=&e
\end{array}\right.$$

\item[c)] Supposons que $r\geq 3$ et que les matrices $A$ et $J$ sont données par 
$$A=
\left(\begin{array}{ccccccccc}
\mu_1&\nu_1&&&&&&&k_1\\
\nu_1&-\mu_1&&&&&&&l_1\\
&&\ddots&&&&&&\vdots\\
&&&\mu_s&\nu_s&&&&k_s\\
&&&\nu_s&-\mu_s&&&&l_s\\
&&&&&\lambda_1&&&h_1\\
&&&&&&\ddots&&\vdots\\
&&&&&&&\lambda_{r-1}&h_{r-1}\\
k_1&l_1&\ldots&k_s&l_s&h_1&\ldots&h_{r-1}&e\\
\end{array}\right)
$$
et 


$$\left(\Id_{1}\bigoplus(-\Id_{1}\right)^s\bigoplus\Id_r$$

Posons alors $J'=(\Id_1\bigoplus(-\Id_1))^{s+1}\bigoplus\Id_{r-2}$ et 

$$A'=
\left(\begin{array}{ccccccccc}
\mu_1&\nu_1&&&&&&&k'_1\\
\nu_1&-\mu_1&&&&&&&l'_1\\
&&\ddots&&&&&&\vdots\\
&&&\mu_{s+1}&\nu_{s+1}&&&&k'_{s+1}\\
&&&\nu_{s+1}&-\mu_{s+1}&&&&l'_{s+1}\\
&&&&&\lambda_3&&&h'_3\\
&&&&&&\ddots&&\vdots\\
&&&&&&&\lambda_{r-1}&h'_{r-1}\\
k'_1&l'_1&\ldots&k'_{s+1}&l'_{s+1}&h'_3&\ldots&h'_{r-1}&e'\\
\end{array}\right)
$$

D'après les formules de la démonstration du Théorème \ref{effective_rep_determ_fleche}, il suffit de prendre les réels $h'_i,k'_j,l'_j$ pour $i=2\ldots r-1$, $j=1\ldots s+1$ et $e'$ tels que :

$$\left\{\begin{array}{rll}
(h'_i)^2&=&h_i^2\frac{(\lambda_i-\lambda_1)(\lambda_i-\lambda_2)}{\left(-(\lambda_i-\mu_{s+1})^2-\nu_{s+1}^2\right)}\\
(k'_j+l'_jI)^2&=&(k_j+l_jI)^2 \frac{(\mu_j+I\nu_j-\lambda_1)(\mu_j+I\nu_j-\lambda_2)}{\left(-(\mu_j+I\nu_j-\mu_{s+1})^2-\nu_{s+1}^2)\right)}\\
(k'_{s+1}-l'_{s+1}I)^2&=&q(\mu_{s+1}+I\nu_{s+1})\times\\
&&\left((I\nu_{s+1})\times\prod_{i=3}^{r-1}(\mu_{s+1}+I\nu_{s+1}-\lambda_i)\times \prod_{j=1}^s(-(\mu_{s+1}-\mu_j+I\nu_{s+1})^2-\nu_j^2)\right)^{-1}\\
e'&=&e+2\mu_{s+1}+4\sum_{j=1}^s\mu_j-2\sum_{i=3}^{r-1}-\lambda_1-\lambda_2\\
\end{array}\right.$$


\ie
Ce qui termine la démonstration.
\end{preuve}


Ainsi, on obtient une sorte de hiérarchie des représentations déterminantales suivant leur signature. Sachant que, d'une part tout polynôme de degré $d$ possède une représentation de signature $(\lfloor\frac{d}{2}\rfloor,d-\lfloor\frac{d}{2}\rfloor)$, et d'autre part pour obtenir une représentation déterminantale de signature plus "grande" il faut disposer d'une minoration du nombre des racines réelles.\par
Pour poursuivre dans la même veine, on peut bien entendu s'interroger sur la manière de passer formellement d'une représentation déterminantale (pas forcément flèche diagonale) à un autre de signature inférieure, et aussi sur l'existence d'une telle hiérarchie dans le cas de plusieurs variables. 



\end{document}